\documentclass{amsart}
\usepackage{latexsym, amssymb, amsmath, amsthm, rotating, amscd}
\usepackage{graphics}
\newtheorem{theorem}{Theorem}[section]
\newtheorem{definition}[theorem]{Definition}

\newtheorem{lemma}[theorem]{Lemma}
\newtheorem{proposition}[theorem]{Proposition}

\numberwithin{equation}{section}

\def\Z{\mathbb{Z}}
\def\N{\mathbb{N}}
\def\C{\mathbb{F}}

\def\co{\begin{sideways}\begin{sideways}$Y$\end{sideways}\end{sideways}}
\def\no{ _{\circ}^{\circ}}

\begin{document}

\title[Constructions of VOCs]{Constructions of vertex operator 
coalgebras via vertex operator algebras}
\author{Keith Hubbard}
\date{\today}

\begin{abstract}
The notion of vertex operator coalgebra is presented which corresponds to the
family of correlation functions of one string propagating in space-time
splitting into n strings in conformal field theory.  This notion is in some
sense dual to the notion of vertex operator algebra.  We prove that any vertex
operator algebra equipped with a non-degenerate, Virasoro preserving, bilinear
form gives rise to a corresponding vertex operator coalgebra.
\end{abstract}

\maketitle

\section{Introduction}
The theory of vertex operator algebras has been an ever expanding field since its inception 
in the 1980s when Borcherds first introduced the precise notion of a vertex algebra (\cite{B}).  
Specialized to vertex 
operator algebras (VOAs) by Frenkel, Lepowsky, and Meurman in \cite{FLM} and shown to have 
inherent connections to conformal field theory,
modular forms, finite groups and Lie theory, vertex (operator) algebras have gained
interest throughout the mathematical community.  Subsequently VOAs have been interpreted 
as a specific type of algebra over an operad (\cite{HL}), and have possible applications 
in giving a geometric definition of elliptic cohomology (\cite{ST}).  Quite recently, in
\cite{K} and \cite{K1}, the algebra induced by considering the operad structure of
worldsheets swept out by closed strings propagating through space-time, has been supplemented by 
examining the induced coalgebra structure which gives rise to vertex operator coalgebras (VOCs).

The notion of VOC corresponds to the coalgebra of correlation functions of one string splitting
into $n$ strings in space-time, whereas VOAs correspond to the algebra of correlation functions 
of $n$ strings combining into one string in space-time.  Not only are the notions of VOC and
VOA integral parts of conformal field theory but the successful integration of the two notions 
would provide an understanding of the correlation functions
arising from any 2-dimensional worldsheet, including higher genus worldsheets.

Examples of VOAs have been quite useful, both in understanding VOA structure and also in understanding 
other objects (particularly the Monster finite simple group \cite{FLM}), but have historically been quite 
challenging to construct.  Building a substantial pool of examples has taken decades.  (See \cite{LL} for 
one recent list of constructions.)  One might expect similar challenges in generating examples of VOCs.  
However, in this paper, via an appropriately defined bilinear form, we generate a large family of 
examples of VOCs by tapping into the extensive work on VOA examples.  We will also explicitly calculate
the family of examples corresponding to Heisenberg algebras and demonstrate the naturally adjoint nature
of Heisenberg VOAs and VOCs.

Finally, dating back to the 1980s conformal field theories have included
axioms about the adjointness of operators induced by manifolds of opposite
orientation (cf. \cite{S}).  The construction in this paper indicates that VOC
operators satisfy this condition and are adjoint to VOA operators when an
appropriate bilinear form exists

The author is indebted to his advisor, Katrina Barron, for her guidance and encouragement throughout 
the research, writing and revising phases of this paper which corresponds to the final chapter of his 
Ph.D. thesis.  
In addition, he thanks Stephan Stolz for numerous helpful conversations.  
Discussions with Bill Dwyer, Brian Hall, 
Hai-sheng Li and Corbett Redden were also much appreciated.  The author gratefully acknowledges the 
financial support of the Arthur J. Schmidt Foundation.

\section{Definitions and algebraic preliminaries}

We begin by reviewing a necessary series from the calculus of formal variables,
then recall the definition of vertex operator coalgebra (cf. \cite{K}, \cite{K1}) 
and vertex operator algebra (cf. \cite{FLM}, \cite{FHL}).  For later use, we will include two standard 
consequences of the definition of VOAs.

\subsection{Delta functions}

We define the ``formal $\delta$-function" to be
\begin{equation*}
\delta(x)= \sum_{n \in \Z} x^n.
\end{equation*}

\noindent
Given commuting formal variables $x_1$, $x_2$ and $n$ an integer, $(x_1 \pm x_2)^n$ will be understood to be
expanded in nonnegative integral powers of $x_2$.  (This is the convention throughout vertex operator algebra
and coalgebra theory.)  Note that the $\delta$-function applied to 
$\frac{x_1-x_2}{x_0}$, where $x_0$, $x_1$ and $x_2$ are commuting formal 
variables, is a formal power series in nonnegative integral powers of $x_2$ 
(cf. \cite{FLM}, \cite{FHL}).

\subsection{The notion of vertex operator coalgebra}
The following description of a vertex operator coalgebra is the central structure of this paper.
Originally motivated by the geometry of propagating strings in conformal field theory 
(\cite{K}, \cite{K1}), VOCs may be formulated in terms of vectors spaces over an 
arbitrary field $\C$ using formal commuting variables $x$, $x_1$, $x_2$, $x_3$.

\begin{definition}\label{D:voc}
A \emph{vertex operator coalgebra (over $\C$) of rank $d \in \C$} is
a $\Z$-graded vector space over $\C$

\begin{equation*}
V = \coprod_{k \in \Z} V_{(k)}
\end{equation*} 

\noindent
such that $\dim V_{(k)} < \infty$ for $k \in \Z$ and
$V_{(k)} = 0$ for $k$ sufficiently small,
together with linear maps

\begin{align*}
\co (x) : V &\mapsto (V \otimes V)[[x,x^{-1}]] \\
v &\mapsto \co(x)v = \sum_{k\in \Z} \Delta_k(v) x^{-k-1}
\ \ \ \ \ (\text{where }\Delta_k(v) \in V \otimes V),
\end{align*}

\begin{equation*}
 c : V \mapsto \C,
\end{equation*}

\begin{equation*}
\rho : V \mapsto \C,
\end{equation*}

\noindent
called the \emph{coproduct}, the \emph{covacuum map} and the \emph{co-Virasoro map}, respectively,
satisfying the following 7 axioms:

1. Left Counit: For all $v \in V$ 

\begin{equation}\label{E:counit}
(c \otimes Id_V) \co(x)v=v
\end{equation} 

2. Cocreation: For all $v \in V$

\begin{equation} \label{E:cocreat1}
(Id_V \otimes c) \co(x)v \in V[[x]] \ \ \text{and}
\end{equation}
\begin{equation} \label{E:cocreat2}
\lim_{x \to 0} (Id_V \otimes c) \co(x)v=v.
\end{equation}

3. Truncation: Given $v \in V$, then $\Delta_k(v) = 0$ for $k$ sufficiently small.

4. Jacobi Identity:

\begin{multline} \label{E:Jac}
x_0^{-1}\delta \left(\frac{x_1-x_2}{x_0} \right) 
(Id_V \otimes \co(x_2))  \co(x_1)
-x_0^{-1}\delta \left(\frac{x_2-x_1}{-x_0} \right)
(T \otimes Id_V)\\
 (Id_V \otimes \co(x_1))  \co(x_2) 
=x_2^{-1}\delta \left(\frac{x_1-x_0}{x_2} \right)
(\co(x_0) \otimes Id_V)  \co(x_2).
\end{multline}

5. Virasoro Algebra:
The Virasoro algebra bracket,

\begin{equation*}
[L(j),L(k)]=(j-k)L(j+k)+\frac{1}{12}(j^3-j)\delta_{j,-k}d,
\end{equation*}

\noindent
holds for $j, k \in \Z$, where 

\begin{equation} \label{E:L_def}
(\rho \otimes Id_V) \co(x) = \sum_{k \in \Z} L(k) x^{k-2}.
\end{equation}

6. Grading:  For each $k \in \Z$ and $v \in V_{(k)}$ 

\begin{equation} \label{E:VOCgrading}
L(0)v= kv.
\end{equation}

7. $L(1)$-Derivative:

\begin{equation} \label{E:L(1)deriv}
\frac{d}{dx} \co(x)=
(L(1) \otimes Id_V) \co(x).
\end{equation}

\end{definition}

\noindent
We denote this vertex operator coalgebra by 
$(V, \co, c, \rho)$ or simply by $V$ when the structure is clear.

Note that $\co$ is linear so that, for example,
$(Id_V \otimes \co(x_1))$ acting on the coefficients of $\co(x_2)v \in
(V \otimes V)[[x_2,x_2^{-1}]]$ is well defined. Notice also, that when each 
expression is applied to any element of $V$, the coefficient of each 
monomial in the formal variables is a finite sum.

\subsection{The definition of a VOA}  \label{S:VOAdef}

Vertex operator algebras were first defined in \cite{FLM}, but it was not until
\cite{H2}, \cite{H} that this definition was rigorously tied to the geometry of conformal
field theory.  This correspondence, along with its operadic interpretation in \cite{HL}, 
helped to motivate the notion of VOC in \cite{K}.  It is not surprising then, that the 
axioms of VOA and VOC would strongly resemble each other.

\begin{definition}\label{D:voa}
A \emph{vertex operator algebra (over $\C$) of rank $d \in \C$} is
a $\Z$-graded vector space over $\C$

\begin{equation*}
V = \coprod_{k \in \Z} V_{(k)}, 
\end{equation*} 

\noindent
such that $\dim V_{(k)} < \infty$ for $k \in \Z$ and
$V_{(k)} = 0$ for $k$ sufficiently small,
together with a linear map $V \otimes V \to V[[x,x^{-1}]]$, or equivalently,

\begin{align*}
Y (\cdot,x): V &\mapsto (\text{End } V) [[x,x^{-1}]] \\
v &\mapsto Y(v,x) = \sum_{k\in \Z} v_k x^{-k-1} 
\ \ \ \ \text{ (where } v_n \in \text{End }V),
\end{align*}

\noindent
and equipped with two distinguished homogeneous vectors in $V$, $\mathbf{1}$ 
(the \emph{vacuum}) and $\omega$ (the \emph{Virasoro element}),
satisfying the following 7 axioms:

1. Left unit: For all $v \in V$ 

\begin{equation}\label{E:unit}
Y(\mathbf{1},x)v=v
\end{equation} 

2. Creation: For all $v \in V$

\begin{equation} \label{E:creat1}
Y(v,x) \mathbf{1} \in V[[x]] \ \ \text{and}
\end{equation}
\begin{equation} \label{E:creat2}
\lim_{x \to 0} Y(v,x) \mathbf{1} =v.
\end{equation}

3. Truncation: Given $v,w \in V$, then $v_k w = 0$ for $k$ sufficiently large.

4. Jacobi Identity: For all $u,v \in V$,

\begin{multline} \label{E:Jac2}
x_0^{-1}\delta \left(\frac{x_1-x_2}{x_0} \right) 
Y(u,x_1)Y(v,x_2)
-x_0^{-1}\delta \left(\frac{x_2-x_1}{-x_0} \right)
Y(v,x_2)Y(u,x_1) \\
=x_2^{-1}\delta \left(\frac{x_1-x_0}{x_2} \right)
Y(Y(u,x_0)v,x_2).
\end{multline}

5. Virasoro Algebra:
The Virasoro algebra bracket,

\begin{equation*}
[L(j),L(k)]=(j-k)L(j+k)+\frac{1}{12}(j^3-j)\delta_{j,-k}d,
\end{equation*}

\noindent
holds for $j, k \in \Z$, where 

\begin{equation} \label{E:L_def1}
Y(\omega,x) = \sum_{k \in \Z} L(k) x^{-k-2}.
\end{equation}

6. Grading:  For each $k \in \Z$ and $v \in V_{(k)}$ 

\begin{equation} \label{E:VOCgrading1}
L(0)v= kv.
\end{equation}

7. $L(-1)$-Derivative: Given $v \in V$,

\begin{equation} \label{E:L(-1)deriv}
\frac{d}{dx} Y(v,x)=
Y(L(-1)v,x).
\end{equation}

\end{definition}

We denote a VOA either by $V$ or by the quadruple $(V, Y, \mathbf{1}, \omega)$.
A vector $v \in V_{(k)}$ for some $k \in \Z$ is said to be a \emph{homogeneous vector
of weight $k$} and we write $\text{wt }v=k$.
A pair of basic properties of VOAs will be necessary for our discussion (cf. 
\cite{FHL}, \cite{LL}):

\begin{align*} 
Y(v,x) \mathbf{1} &=e^{xL(-1)}v \ \ \ \ \ \ \ \ \ \ \ \ \ \  \text{ for } v \in V, \\
\text{wt } v_k w&=r + s-k-1     \ \ \ \ \ \ \ \  \text{ for } v \in V_{(r)}, \ w \in V_{(s)}
\end{align*}

\section{A family of examples of VOCs}
One of the most natural questions to ask about vertex operator coalgebras is
``what do they
look like?",  or even, ``do any exist?".  The main purpose of this paper is
to answer the latter question in the affirmative and to provide concrete insight
into the former question.

\subsection{A family of examples of VOCs}  \label{S:family_examples}
Let the vector space $V = \coprod_{k \in \Z} V_{(k)}$ be a module over the Virasoro
algebra, $\mathcal{V} = \oplus_{j \in \Z} \C L(j) \oplus \C d$, such that for 
all homogeneous vectors $L(0) \cdot v= \text{wt }(v) v$.  
We will say that a bilinear form $(\cdot,\cdot)$ on $V$ is \emph{Virasoro preserving} 
if it satisfies the condition

\begin{equation} \label{E:L-antisym}
(L(k)v_1,v_2)=(v_1,L(-k)v_2)
\end{equation}

\noindent
for all $k \in \Z$, $v_1,v_2 \in V$.  In particular, $k=0$ in Property (\ref{E:L-antisym})
indicates that all Virasoro preserving bilinear forms are graded, i.e.

\begin{equation*}
(V_{(k)},V_{(\ell)})=0
\end{equation*}

\noindent
for $k \neq \ell$.  If $V$ is a VOA, the bilinear form is said to be  
\emph{invariant} if, for all $u,v,w \in V$,

\begin{equation} \label{E:VOAinvar}
(Y(u,x)v,w)=(u,Y(e^{xL(1)}(-x^{-2})^{L(0)}v,x^{-1})w).
\end{equation}

\noindent
Any invariant bilinear form on V is 
Virasoro preserving ((2.31) in \cite{L}).

Note that there is a natural extension of $(\cdot,\cdot):V^{\otimes 2} \to \C$ to
$(\cdot,\cdot):V^{\otimes 4} \to \C$ given by 
$(u_1 \otimes u_2,v_1 \otimes v_2)=(u_1,v_1) (u_2,v_2)$, for $u_1,u_2,v_1,v_2 \in V$.

\begin{theorem} \label{T:examples}
Let $(V,Y,\mathbf{1},\omega)$ be a vertex operator algebra equipped with a
nondegenerate and Virasoro preserving bilinear form $(\cdot,\cdot)$.
Given the linear operators

\begin{align*}
 c : V & \to \C \\
v &\mapsto (v,\mathbf{1}),
\end{align*}

\begin{align*}
 \rho : V & \to \C \\
v &\mapsto (v,\omega),
\end{align*}

\noindent
and 

\begin{align*}
\co : V &\to (V \otimes V)[[x,x^{-1}]] \\
v &\mapsto \co(x)v = \sum_{k\in \Z} \Delta_k(v) x^{-k-1},
\end{align*}

\noindent
defined by

\begin{equation}
(\co(x)u, v \otimes w) = (u, Y(v, x)w),
\end{equation}

\noindent 
the quadruple $(V,\co,c, \rho)$ is a vertex operator coalgebra.
\end{theorem}

\begin{proof}  We will show that all 8 axioms for VOCs are satisfied.

1. Positive energy: Trivially satisfied.

2. Left Counit: Given $u \in V$, for all $v \in V$ 

\begin{align*}
((c \otimes Id_V) \co(x)u,v)&=(\co(x)u, \mathbf{1} \otimes v)\\
&=(u, Y(\mathbf{1},x)v) \\
&=(u,v).
\end{align*} 

\noindent
Thus, by nondegeneracy, $(c \otimes Id_V) \co(x)u=u$.

3. Cocreation: Given $u \in V$, then for all $v \in V$ 

\begin{align*}
((Id_V \otimes c) \co(x)u,v)&= (\co(x)u, v \otimes \mathbf{1}) \\
&= (u, Y(v,x)\mathbf{1}) \\
&= (u, e^{x L(-1)}v) \in \C[x]
\end{align*} 

\noindent
and 

\begin{equation*}
\lim_{x \to 0} ( u, e^{x L(-1)}v)= (u,v).
\end{equation*}

4. Truncation: Pick $N \in \Z$ such that $V_{(n)}=0$ for all $n \leq N$.
Given $u \in V_{(r)}$, let $v \in V_{(s)}$ and $w \in V_{(t)}$.

\begin{align*}
( \co(x)u,v \otimes w)
&= ( u,Y(v,x) w)\\
&= \sum_{k \in \Z} (u,v_k w) x^{-k-1}
\end{align*}

\noindent
For $(u,(v)_k w)$ to be nonzero, we must have $\text{wt } u=
\text{wt } v + \text{wt } w -k -1$, i.e. $r=s+t-k-1$;
but $s,t > N$ so we must have $r > 2N-k-1$ or $r-2N > -k-1$.
Hence, $( \co(x)u,v \otimes w) \in \C[[x^{-1}]]x^{r-2N-1}$
for any $s,t \in \Z$.

5. Jacobi Identity: Given $u \in V$, then for all $v_1,v_2,v_3 \in V$ 

\begin{align} \label{E:511}
((Id_V \otimes \co(x_2))  \co(x_1)u,v_1 \otimes v_2 \otimes v_3)
&= (\co(x_1)u,v_1 \otimes Y(v_2,x_2) v_3) \\
&= (u,Y(v_1,x_1) Y(v_2,x_2) v_3), \notag
\end{align}

\begin{align} \label{E:512}
((T \otimes Id_V)
 (Id_V \otimes \co(x_1))  \co(x_2)u, & v_1 \otimes v_2 \otimes v_3) \\
&= ((Id_V \otimes \co(x_1))  \co(x_2)u,v_2 \otimes v_1 \otimes v_3) \notag \\
&= (\co(x_2)u,v_2 \otimes Y(v_1,x_1) v_3) \notag \\
&= (u,Y(v_2,x_2) Y(v_1,x_1) v_3), \notag
\end{align}

\begin{align} \label{E:513}
((\co(x_0) \otimes Id_V)  \co(x_2)u,v_1 \otimes v_2 \otimes v_3)
&= ( \co(x_2)u, Y(v_1,x_0) v_2 \otimes v_3) \\
&= (u,Y(Y(v_1,x_0)v_2,x_2)v_3). \notag
\end{align}

Equations (\ref{E:511}), (\ref{E:512}) and (\ref{E:513}) make it clear that
the VOA Jacobi identity (\ref{E:Jac2}) is equivalent to the VOC Jacobi identity
(\ref{E:Jac}).

6. Virasoro Algebra: Given $u \in V$, for all $v \in V$ 

\begin{align} \label{E:519}
((\rho \otimes Id_V) \co(x)u,v)
&=(\co(x)u, \omega \otimes v) \\
&=(u, Y(\omega,x)v) \notag \\
&=\sum_{k \in \Z} (u, L(k)v) x^{-k-2} \notag \\
&=\sum_{j \in \Z} (L(j)u,v) x^{j-2}. \notag
\end{align}

\noindent
Note that in the last equality we have used Virasoro preservation, 
(\ref{E:L-antisym}).

Equation (\ref{E:519}) shows that the Virasoro algebra bracket,

\begin{equation*}
[L(j),L(k)]=(j-k)L(j+k)+\frac{1}{12}(j^3-j)\delta_{j,-k}d,
\end{equation*}

\noindent
follows from the Virasoro bracket relation on VOAs. 

7. Grading:  Equation (\ref{E:519}) shows that $L(0)=Res_x x Y(\omega,x)$ so 
grading follows from VOAs.

8. $L(1)$-Derivative: Given $u \in V$, then for all $v,w \in V$ 

\begin{align*}
((L(1) \otimes Id_V) \co(x)u,v \otimes w) 
&=(\co(x)u,L(-1)v \otimes w) \\
&= (u, Y(L(-1)v, x)w) \\
&= \frac{d}{dx}(u, Y(v, x)w) \\
&= \frac{d}{dx}(\co(x)u,v \otimes w) 
\end{align*}

\noindent
Here the first equality uses Virasoro preservation.
\end{proof}

Li showed in \cite{L} that if a simple VOA satisfies 
the condition $L(1)V_{(1)}=0$ then there exists a nondegenerate, invariant bilinear 
form on V.  Thus we are guaranteed a family of VOAs equipped with the type of form
required for Theorem \ref{T:examples}.  Additionally, Heisenberg VOAs may be explicitly equipped with
an appropriate bilinear form and they will be the focus of more concrete discussion in the next 
section.

\subsection{Vertex operator algebras and coalgebras associated with Heisenberg algebras}
While the construction in the last section does describe a family of VOCs, it is not
extremely explicit in nature.  In this section we will explicitly construct VOCs 
from Heisenberg algebras.  We begin with the construction of the vector space for the
Heisenberg VOA following \cite{D}.

Let $\mathbf{h}$ be a finite dimensional vector space equipped with a symmetric, 
nondegenerate bilinear form 
$\langle \cdot, \cdot \rangle$.  Since $\mathbf{h}$ may be considered as an abelian Lie
algebra, let $\hat{\mathbf{h}}$ be the corresponding affine Lie algebra, i.e., let

\begin{equation*}
\hat{\mathbf{h}}=\mathbf{h} \otimes \C [t,t^{-1}] \oplus \C c,
\end{equation*}

\noindent
where $c$ is nonzero, with the Lie bracket defined by

\begin{align*}
[\alpha \otimes t^m,\beta \otimes t^n] &= \langle \alpha, \beta \rangle m \delta_{m,-n}c \\
[\hat{\mathbf{h}},c] &=0 
\end{align*}

\noindent
for $\alpha, \beta \in \mathbf{h}$, $m,n \in \Z$.  There is a natural $\Z$-grading
on $\hat{\mathbf{h}}$ under which $\alpha \otimes t^m$ has weight $-m$ for all
$\alpha \in \mathbf{h}$, and $m \in \Z$, and $c$ has weight 0.  The element $\alpha \otimes t^m$
of $\hat{\mathbf{h}}$ is usually denoted $\alpha(m)$.  Three graded subalgebras 
are of interest:

\begin{equation*} 
\hat{\mathbf{h}}^+=\mathbf{h} \otimes t \C [t],
\end{equation*}

\begin{equation*}
\hat{\mathbf{h}}^-=\mathbf{h} \otimes t^{-1} \C [t^{-1}],
\end{equation*}

\begin{equation*}
\hat{\mathbf{h}}_{\Z}=\hat{\mathbf{h}}^+ \oplus \hat{\mathbf{h}}^- \oplus \C c.
\end{equation*}

\noindent
The subalgebra $\hat{\mathbf{h}}_{\Z}$ is a \emph{Heisenberg algebra}, by which we
mean that its center is one-dimensional and is equal to its commutator subalgebra.
Note that $\hat{\mathbf{h}}^+$ and $\hat{\mathbf{h}}^-$ are abelian, but that
$\hat{\mathbf{h}}_{\Z}$ is necessarily non-abelian.

We consider the induced $\hat{\mathbf{h}}_{\Z}$-module

\begin{equation*}
M(1)=U(\hat{\mathbf{h}}_{\Z}) \otimes_{U(\hat{\mathbf{h}}^+ \oplus \C c)} \C
\end{equation*}

\noindent
where $U$ indicates the universal enveloping algebra and $\C$ is viewed as a 
$\Z$-graded $(\hat{\mathbf{h}}^+ \oplus \C c)$-module by

\begin{align*}
c \cdot 1 &= 1, \\
\hat{\mathbf{h}}^+ \cdot 1 &= 0, \\
\text{deg } 1 &=0.
\end{align*}

\noindent
The module $M(1)$ may be generalized (cf. \cite{L} and \cite{FLM}).  
$M(1)$ is linearly isomorphic to $S(\hat{\mathbf{h}}^-)$ in
a way that preserves grading.  Thus we often write basis elements of $M(1)$ as

\begin{equation*}
v=\alpha_1(-n_1) \cdots \alpha_r(-n_r)
\end{equation*}

\noindent
for $\alpha_i \in \mathbf{h}$, $n_i \in \Z_+$, $i=1, \ldots, r$, and observe that $v$ 
has weight $n_1+ \cdots+n_r$.  (Tensor products are suppressed in this notation.)  
Note that the $\alpha_i(-n_i)$'s all commute so
their order is irrelevant.  Using the form $\langle \cdot, \cdot \rangle$ on 
$\mathbf{h}$, we may choose $\{ \gamma_i \}_{i=1}^{d}$ to be an orthonormal basis and we
lose no generality by considering only $\alpha_i$ from this set of basis elements.
Thus, we will typically prove results for the set of generating elements

\begin{equation*}
\text{gen}M = \left\{ \alpha_1(-n_1) \cdots \alpha_r(-n_r) |
r \in \N, \ \alpha_j \in \{ \gamma_i \}_{i=1}^{d}, \ n_j \in \Z_+, \ j=1, \ldots, r \right\}
\end{equation*}

\noindent
and then extend linearly to all of $M(1)$.  

Next we define a bilinear form on $M(1)$ which we will use throughout the rest of this
section.

\begin{lemma}
There is a unique bilinear form $(\cdot,\cdot)$ on $M(1)$ satisfying 

\begin{align}
(\alpha(m) \cdot u, v)&= (u, \alpha(-m) \cdot v),  \label{E:311} \\
(1,1)&=1. \label{E:3181}
\end{align}

\noindent
for all $u,v \in M(1)$, $\alpha \in \mathbf{h}$, $m \in \Z \smallsetminus \{0 \}$.

More precisely, given $v=\alpha_1(-n_1) \cdots \alpha_r(-n_r)$ let

\begin{equation*}
p(v) = \alpha_r(n_r) \cdots \alpha_1(n_1) \alpha_1(-n_1) \cdots \alpha_r(-n_r) \in \Z_+ \subset M(1).
\end{equation*}

\noindent
The unique bilinear form $(\cdot, \cdot)$ on $M(1)$ satisfying (\ref{E:311}) and 
(\ref{E:3181}) is defined on basis elements 
$u, v \in \text{gen}M$ by

\begin{align} \label{E:form}
(u,v)= \left\{ \begin{array}{ll}
p(u) & \text{if } u=v \\ 
0  & \text{otherwise} 
\end{array} \right. 
\end{align}

Further, this form is nondegenerate, graded and symmetric.
\end{lemma}

\begin{proof}
We will construct the form in (\ref{E:form}) from (\ref{E:311}) and (\ref{E:3181}), thus showing 
the form is unique.  Consider $u,v \in \text{gen}M$ such that $u \neq v$.  Then there is an
element $\alpha(-n) \in \hat{\mathbf{h}}$ and a positive integer $t$ such that 
$\alpha(-n)^t$ is contained in $u$ or $v$ but not in the other.  We may assume $\alpha(-n)^t$ 
is in $u$, say $u=\alpha(-n)^t \alpha_1(-n_1) \cdots \alpha_r(-n_r)$.  But then

\begin{eqnarray} 
(u,v) 
&=& ( \alpha(-n)^t \alpha_1(-n_1) \cdots \alpha_r(-n_r),v) \label{E:3311}\\
&=& ( \alpha_1(-n_1) \cdots \alpha_r(-n_r), \alpha(n)^t v) \nonumber \\
&=& (\alpha_1(-n_1) \cdots \alpha_r(-n_r), 0) \nonumber \\
&=&0 \nonumber .
\end{eqnarray} 

\noindent
(As a biproduct, (\ref{E:3311}) shows the form is graded and symmetric.)
Now we need only examine the form applied to a single basis element.
Let $u=\alpha_1(-n_1) \cdots \alpha_r(-n_r)$.  Then

\begin{eqnarray} 
(u,u) 
&=&(\alpha_1(-n_1) \cdots \alpha_r(-n_r), \alpha_1(-n_1) \cdots \alpha_r(-n_r))  \label{E:3312}\\
&=&(\alpha_r(n_r) \cdots \alpha_1(n_1) \alpha_1(-n_1) \cdots \alpha_r(-n_r), 1) \nonumber \\
&=& (p(u),1) \nonumber \\
&=& p(u), \nonumber 
\end{eqnarray} 

\noindent
thus proving that (\ref{E:form}) is the unique form satisfying (\ref{E:311}) and (\ref{E:3181}).  
Nondegeneracy is immediate from (\ref{E:3312}).
\end{proof}

Given $\alpha \in \mathbf{h}$, we will need the 
following three series in formal variable $x$ with coefficients in $\hat{\mathbf{h}}_{\Z}$:

\begin{align*}
\alpha^+(x)   &= \sum_{k \in \Z_+} \alpha(k) x^{-k-1} \\
\alpha^-(x)   &= \sum_{k \in \Z_+} \alpha(-k) x^{k-1} \\
\alpha(x)     &= \alpha^-(x) + \alpha^+(x).
\end{align*} 

\noindent
Via Equation (\ref{E:311}) we see that

\begin{align}
(\alpha^-(x)  v_1, v_2)&= (v_1, \alpha^+(x)  v_2),  \label{E:313} \\
(\alpha(x) v_1, v_2)&= (v_1, \alpha(x) v_2).  \label{E:315}
\end{align}

The vertex operator algebra associated to a Heisenberg algebra is
described using a normal ordering procedure, indicated by open colons 
$\no \ \ \no$,
which reorders the enclosed expression so that all operators
$\alpha_1(-m)$ are placed to the left (multiplicatively) of all 
operators $\alpha_2(n)$, for $\alpha_1, \alpha_2 \in \mathbf{h}$, $m, n \in \Z_+$.  For example,

\begin{align*}
\no \alpha_1(x) \alpha_2(x)  \no v
&= \ \no  (\alpha_1^-(x) + \alpha_1^+(x))
   \alpha_2(x)  \no v \\ 
&=   \alpha_1^-(x) \alpha_2(x) v + 
   \alpha_2(x)  \alpha_1^+(x) v.
\end{align*}

\noindent
This normal ordering moves degree-lowering operators to 
the right.  Notice that $\alpha^+(x)$ will only produce elements of lower weight 
than the basis element of $M(1)$ to which it is applied,
while $\alpha^-(x)$ will only produce elements of higher weight.  Therefore 
applying all the degree-lowering operators before all the degree-raising operators 
guarantees that no single 
weight-space has infinitely many summands in it.  

We now have the notation to describe the VOA associated to a Heisenberg algebra.  
Define a linear map $Y ( \cdot,x): M(1) \to (\text{End } M(1))[[x,x^{-1}]]$ by

\begin{equation*}
Y(v, x) = \ \no
\left( \frac{1}{(n_1-1)!} \left(\frac{d}{dz}\right)^{n_1-1} \alpha_1(x) \right) \cdots
\left( \frac{1}{(n_r-1)!} \left(\frac{d}{dz}\right)^{n_r-1} \alpha_r(x) \right)
\ \no 
\end{equation*}

\noindent
for $v=\alpha_1(-n_1) \cdots \alpha_r(-n_r)$.  We also define two 
distinguished elements of $M(1)$, $\mathbf{1}=1$ and 
$\omega = \frac{1}{2} \sum_{i=1}^{d} \gamma_i(-1)^2$, where $\{ \gamma_i \}_{i=1}^{d}$
is the orthonormal basis of $\mathbf{h}$ as above. 

\begin{proposition}
The quadruple $(M(1), Y, \mathbf{1}, \omega)$ as defined above is a vertex operator
algebra.
\end{proposition}

Our treatment here largely mirrors \cite{D}, with the above proposition being 
Proposition 3.1 in \cite{D}.  A proof may be found in \cite{G}. 
For our purposes, however, the nondegenerate bilinear form on $M(1)$, 
described in (\ref{E:form}), is equally relevant.  We will now prove an additional  
property of that bilinear form.

\begin{lemma}
The bilinear form on $M(1)$ defined in Equation (\ref{E:form}) is Virasoro preserving.
\end{lemma}

\begin{proof}
First, we will explicitly calculate the $L(k)$ operators and then show
that $(L(k)v,w)=(v,L(-k)w)$.  Symmetry of the form allows us to only consider $k \in \N$.

Using the definition of the $L(k)$ operators we see that 
$\sum_{k \in \Z} L(k) x^{-k-2} = Y(\frac{1}{2}\sum_{i=1}^{d} \gamma_i(-1)^2, x)$.
Employing the definition of $Y$, for $k \in \Z_+$ we have

\begin{align*}
L(k) &= \sum_{i=1}^{d} \left( \frac{1}{2} \sum_{j =1}^{k-1} \gamma_i(j) \gamma_i(k-j)
+ \sum_{j \in \Z_+} \gamma_i(-j) \gamma_i(k+j) \right) \\
L(-k) &= \sum_{i=1}^{d} \left( \frac{1}{2} \sum_{j =1}^{k-1} \gamma_i(-j) \gamma_i(-k+j)
+ \sum_{j \in \Z_+} \gamma_i(-k-j) \gamma_i(j) \right).
\end{align*}

\noindent
Given $u, v \in \text{gen}M$, 

\begin{equation*}
(L(0)u, v)=\text{wt }(u) \ (u,v) = \text{wt }(v) \ (u,v) = (u,L(0)v)
\end{equation*}

\noindent
since either $(u,v)=0$, or $u=v$ implying $\text{wt }(u)=\text{wt }(v)$.  For $k \in \Z_+$
we make use of (\ref{E:311}) to observe that,

\begin{align*}
(L(k)u,v) 
&=  \sum_{i=1}^{d} \left( \frac{1}{2} \sum_{j =1}^{k-1} (\gamma_i(j) \gamma_i(k-j)u,v)
+  \sum_{j \in \Z_+} (\gamma_i(-j) \gamma_i(k+j)u,v) \right) \\
&= \sum_{i=1}^{d} \left( \frac{1}{2} \sum_{j =1}^{k-1} (u,\gamma_i(-k+j) \gamma_i(-j)v)
+ \sum_{j \in \Z_+} (u,\gamma_i(-k-j) \gamma_i(j)v) \right) \\
&= (u,L(-k)v).
\end{align*}

\end{proof}

Using the Heisenberg VOA $M(1)$ along with the nondegenerate, Virasoro preserving bilinear
form defined in (\ref{E:form}), 
we will follow the construction of a VOC in Section \ref{S:family_examples}.
First, we define a linear map $c : V \to \C$ by

\begin{equation*}
c(1)  =1 
\end{equation*}
\begin{equation*}
c(\alpha_1(-n_1) \cdots \alpha_r(-n_r)) = 0
\end{equation*}

\noindent
where $r \geq 1$ and $\alpha_1(-n_1) \cdots \alpha_r(-n_r) \in \text{gen}M$.  
It is clear that $c(v)=(v,\mathbf{1})$ for all $v \in M(1)$.  
Next, we define a linear map $\rho : V  \to \C$ by

\begin{align*}
\rho(\gamma_i(-1)^2)= 1
\end{align*}

\noindent
for each basis element $\gamma_i$ of $\mathbf{h}$ and $\rho(v)=0$ for $v$ any
other basis element of $M(1)$.  Again it is clear that $\rho(v)=(v,\omega)$ for 
all $v \in M(1)$.  Finally, we need to define a linear map 
$\co (x): V \to (V \otimes V)[[x,x^{-1}]]$ such that

\begin{equation*}
(\co(x)u, v \otimes w) = (u, Y(v, x)w).
\end{equation*}

\noindent
For notational simplicity, given $\alpha \in \mathbf{h}$, $n \in \Z_+$ and $x$ 
a formal variable, let

\begin{align*}
\alpha^+(n,x)   &= \frac{1}{(n-1)!}  \left(\frac{d}{dz}\right)^{n-1}
\sum_{k \in \Z_+} \alpha(k) x^{-k-1} \\
\alpha^-(n,x)   &= \frac{1}{(n-1)!}  \left(\frac{d}{dz}\right)^{n-1}
\sum_{k \in \Z_+} \alpha(-k) x^{k-1} \\
\alpha(n,x)     &= \alpha^-(n,x) + \alpha^+(n,x).
\end{align*}

\begin{proposition}
Let $v$ denote the basis element $\beta_1(-m_1) \cdots \beta_s(-m_s)$ and define
$\co (x): V \to (V \otimes V)[[x,x^{-1}]]$ as

\begin{equation*}
\co(x)u=\sum_{v \in \text{gen}M}
\frac{1}{p(v)} v \otimes \ \no 
\beta_1(m_1,x) \cdots  \beta_s(m_s,x) \ \no u.
\end{equation*}

\noindent
For all $u,v,w \in M(1)$, $(\co(x)u, v \otimes w) = (u, Y(v, x)w)$.
\end{proposition}

\begin{proof}
First, note that our definition of $\co$ is equivalent to 

\begin{equation*}
(\co(x)u, v \otimes w) = 
(\frac{1}{p(v)} v \otimes \ \no \beta_1(m_1,x) \cdots  \beta_s(m_s,x) \ \no u, v \otimes w)
\end{equation*}

\noindent
for all 
$u = \alpha_1(-n_1) \cdots \alpha_r(-n_r)$,
$v = \beta_1(-m_1) \cdots \beta_s(-m_s)$,
$w = \mu_1(-\ell_1) \cdots \mu_t(-\ell_t)$.  Given these basis elements we use 
induction on $s$ to show that

\begin{equation} \label{E:34}
(u, \ \no \beta_1(m_1,x) \cdots  \beta_s(m_s,x) \ \no w) 
=(\ \no \beta_1(m_1,x) \cdots  \beta_s(m_s,x) \ \no u, w).
\end{equation}

\noindent
For $s=0$, this is trivial.  If we assume that (\ref{E:34}) is true for $s-1$ and 
appeal to (\ref{E:313}), we see that

\begin{equation*}
(u, \ \no \beta_1(m_1,x) \cdots  \beta_s(m_s,x) \ \no w)
\end{equation*}
\begin{multline*}
= (u,  \beta_s^-(m_s,x) \ \no \beta_1(m_1,x) \cdots  \beta_s(m_{s-1},x) \ \no w) \\
+ (u, \ \no \beta_1(m_1,x) \cdots  \beta_s(m_{s-1},x) \ \no \ \beta_s^+(m_s,x)  w) 
\end{multline*}
\begin{multline*}
=  ( \beta_s^+(m_s,x) u, \ \no \beta_1(m_1,x) \cdots  \beta_s(m_{s-1},x) \ \no w) \\
+  (\ \no \beta_1(m_1,x) \cdots  \beta_s(m_{s-1},x)\ \no u, \beta_s^+(m_s,x) w)
\end{multline*}
\begin{multline*}
=  ( \ \no \beta_1(m_1,x) \cdots  \beta_s(m_{s-1},x) \ \no \beta_s^+(m_s,x)  u, w) \\
+ \beta_s^-(m_s,x) \ \no \beta_1(m_1,x) \cdots  \beta_s(m_{s-1},x) \ \no u, w) 
\end{multline*}
\begin{equation*}
=(\ \no \beta_1(m_1,x) \cdots  \beta_s(m_s,x) \ \no u, w).
\end{equation*}

\noindent
Finally, using Equation (\ref{E:34}) we see that

\begin{align*}
(\co(x)u, v \otimes w) 
&=(\frac{1}{p(v)} v \otimes \ \no \beta_1(m_1,x) \cdots  \beta_s(m_s,x)
\ \no u, v \otimes w) \\
&=\frac{(v, v)}{p(v)}  ( \ \no \beta_1(m_1,x) \cdots  \beta_s(m_s,x) \ \no u, w) \\
&=(u, \ \no \beta_1(m_1,x) \cdots  \beta_s(m_s,x) \ \no w) \\
&=(u, Y(v, x)w).
\end{align*}
\end{proof}

Theorem \ref{T:examples} proves that the quadruple $(M(1),\co,c,\rho)$ associated 
to the Heisenberg algebra $\hat{\mathbf{h}}_{\Z}$ is a vertex operator coalgebra.

\end{document}